\documentclass{article}

\usepackage{fullpage}

\DeclareMathAlphabet{\mathcal}{OMS}{cmsy}{m}{n}

\makeatletter
\xdef\@endgadget#1{{\unskip\nobreak\hfil\penalty50\hskip1em\hbox{}\nobreak\hfil#1\parfillskip=0pt\finalhyphendemerits=0\par}}
\newcommand\@Endofsymbol{$\triangledown$}
\newcommand\Endofremark{\@endgadget{\@Endofsymbol}}
\makeatother

\title{Tractable downfall of basis pursuit in structured sparse optimization \thanks{The project was supported by the Israel Science Foundation (grant no. 2406/22), while the second author was a Jane and Larry Sherman Fellow.}}

\author{Maya V. Marmary\thanks{M. V. Marmary is with the Faculty of Mechanical Engineering, Technion -- Israel Institute of Technology, Haifa, Israel
		{\tt maya.marmary@campus.technion.ac.il}} \quad and \quad Christian Grussler\thanks{C. Grussler is with the Stephen B. Klein Faculty of Aerospace Engineering, Technion -- Israel Institute of Technology, Haifa, Israel
		{\tt cgrussler@technion.ac.il}} 
}

\usepackage{amsthm}
\usepackage{graphicx}
\usepackage{enumerate}
\usepackage{mathtools}
\usepackage{amsmath} %
\usepackage{xcolor} %

\usepackage{enumitem}

\usepackage{amssymb}  %
\usepackage{bm}

\usepackage{bbm}

\usepackage{tikz}
\usepackage{tkz-euclide}

\definecolor{ao(english)}{rgb}{0.0, 0.5, 0.0}
\usepackage[colorlinks, citecolor = {ao(english)}, linkcolor = {ao(english)}, pdfkeywords={Basis pursuit, total positivity, sparse optimization}]{hyperref} 
\usepackage[capitalize,nameinlink]{cleveref}
\usepackage{siunitx}
\usepackage{pgfplots}
\pgfplotsset{compat=newest} 
\pgfplotsset{compat=newest}

\usetikzlibrary{arrows,shapes,trees,calc,positioning,patterns,decorations.pathmorphing,decorations.markings,backgrounds}
\usetikzlibrary{matrix}

\usepgfplotslibrary{groupplots}

\newtheorem{thm}{Theorem}
\crefname{thm}{Theorem}{Theorems}

\newtheorem{prop}{Proposition}
\crefname{prop}{Proposition}{Propositions}

\newtheorem{lem}{Lemma}
\crefname{lem}{Lemma}{Lemmas}

\newtheorem{cor}{Corollary}
\crefname{cor}{Corollary}{Corollaries}

\newtheorem{rem}{Remark}
\crefname{rem}{Remark}{Remark}

\crefname{ass}{Assumption}{Assumption}

\crefname{conj}{Conjecture}{Conjectures}

\newtheorem{defn}{Definition}
\crefname{defn}{Definition}{Definitions}

\crefname{prob}{Problem}{Problems}

\crefname{appl}{Application}{Applications}

\newtheorem{exm}{Example}
\crefname{exm}{Example}{Examples}

\crefname{algorithm}{Algorithm}{Algorithms}
\crefname{paper}{Paper}{Papers}
\crefname{figure}{Figure}{Figures}
\crefname{section}{Section}{Sections}
\Crefname{section}{Section}{Sections}

\usepackage{dsfont}
\let\mathbb=\mathds

\newcommand{\Rmn}{\mathbb{R}^{m \times n}}

\DeclareMathOperator{\sign}{sign} %
\newcommand{\rk}{\textnormal{rank}}

\newcommand{\diag}{\textnormal{diag}}

\newcommand{\vari}[1]{\text{S}(#1)}

\newcommand{\vb}[1]{\text{VB}_{#1}}

\newcommand{\sgc}[1]{{SC}_{#1}}

\newcommand{\normA}[2]{\| {} #1{} \|_{ #2}}

\newcommand{\transp}{\mathsf{T}}

\colorlet{FigColor1}{blue}
\colorlet{FigColor2}{red}
\colorlet{FigColor3}{ao(english)}
\colorlet{FigColor4}{orange}
\pgfplotsset{every axis plot/.append style={line width=1.5pt}}
\definecolor{bluebell}{rgb}{0.74, 0.83, 0.9}
\definecolor{airforceblue}{rgb}{0.36, 0.54, 0.66}

\crefformat{equation}{\textup{#2(#1)#3}}
\crefrangeformat{equation}{\textup{#3(#1)#4--#5(#2)#6}}
\crefmultiformat{equation}{\textup{#2(#1)#3}}{ and \textup{#2(#1)#3}}
{, \textup{#2(#1)#3}}{, and \textup{#2(#1)#3}}
\crefrangemultiformat{equation}{\textup{#3(#1)#4--#5(#2)#6}}%
{ and \textup{#3(#1)#4--#5(#2)#6}}{, \textup{#3(#1)#4--#5(#2)#6}}{, and \textup{#3(#1)#4--#5(#2)#6}}

\Crefformat{equation}{#2Equation~\textup{(#1)}#3}
\Crefrangeformat{equation}{Equations~\textup{#3(#1)#4--#5(#2)#6}}
\Crefmultiformat{equation}{Equations~\textup{#2(#1)#3}}{ and \textup{#2(#1)#3}}
{, \textup{#2(#1)#3}}{, and \textup{#2(#1)#3}}
\Crefrangemultiformat{equation}{Equations~\textup{#3(#1)#4--#5(#2)#6}}%
{ and \textup{#3(#1)#4--#5(#2)#6}}{, \textup{#3(#1)#4--#5(#2)#6}}{, and \textup{#3(#1)#4--#5(#2)#6}}

\crefdefaultlabelformat{#2\textup{#1}#3}

\usepackage[style=numeric, backend=biber,sortcites=true]{biblatex}
\providecommand{\keywords}[1]{
		
	\textbf{\textit{Keywords---}} #1
}
\begin{document}

\maketitle
 \begin{abstract}
	The problem of finding the sparsest solution to a linear underdetermined system of equations, often appearing, e.g., in data analysis, optimal control, system identification or sensor selection problems, is considered. This non-convex problem is commonly solved by convexification via $\ell_1$-norm minimization, known as basis pursuit (BP). In this work, a class of structured matrices, representing the system of equations, is introduced for which (BP) tractably fails to recover the sparsest solution. In particular, this enables efficient identification of matrix columns corresponding to unrecoverable non-zero entries of the sparsest solution and determination of the uniqueness of such a solution. These deterministic guarantees complement popular probabilistic ones and provide insights into the a priori design of sparse optimization problems. As our matrix structures appear naturally in optimal control problems, we exemplify our findings based on a fuel-optimal control problem for a class of discrete-time linear time-invariant systems. Finally, we draw connections of our results to compressed sensing and common basis functions in geometric modeling. 
\end{abstract}

 
\pagestyle{plain}
\thispagestyle{plain}

\section{Introduction} \label{sec:intro}
{S}{parse} optimization problems are fundamental in various areas, including signal processing \cite{elad2010sparse}, machine learning \cite{sra2011optimization}, compressed sensing \cite{candes2008introduction,donoho2006compressed}, system identification \cite{chen2009sparse,brunton2022data}, optimal control \cite{athans1966optimal,Nagahara2016handsoff}, sparse controller design \cite{fardad2011sparsity} or actuator/sensor selection \cite{zarjovCDC18,zarmohdhigeojovTAC20}.These problems involve finding a feasible solution with the fewest non-zero entries. For instance, in signal processing, sparse representations enable efficient compression by focusing on the most significant components of a signal, while in control applications, sparse solutions simplify models or control laws, leading to more practical, resource-friendly implementations.

A frequently studied setup is the determination of the sparsest solution to an underdetermined system of linear equations \cite[Section~\S 1]{elad2010sparse}, mathematically defined as:
\begin{equation} \label{eq:ell_0}\tag{\text{$\ell_0$}}
	\begin{aligned}
		& \underset{u \in \mathds{R}^n}{\text{minimize}}
		& & \normA{u}{l_0} := \sum_{i=1}^n |\sign(u_i)| \\
		& \text{subject to}
		& & Vu = y
	\end{aligned}
\end{equation}
for given $V\in\mathbb{R}^{m \times n}$, $m < n$, $y \in \mathds{R}^n$ and $\sign(\cdot)$ denoting the sign function taking values in $\{-1,0,1\}$.   
Unfortunately, as $\|\cdot\|_{\ell_0}$ is non-convex, determining the solution to this problem is generally NP-hard \cite{foucart2013invitation}. To this end, convex relaxations have been developed over the past decades \cite{foucart2013invitation,vidyasagar2019introduction,elad2010sparse}. The most common one, known as the Basis Pursuit problem (BP) \cite{chen1994basis}, replaces the $\ell_0$-norm with the $\ell_1$-norm:%
\begin{equation*}\label{eq:ell_1}\tag{\text{BP}}
	\begin{aligned}
		& \underset{u \in \mathds{R}^n}{\text{minimize}}
		& & \normA{u}{l_1} := \sum_{i=1}^n |u_{i}| \\
		& \text{subject to}
		& & Vu = y.
	\end{aligned}
\end{equation*}
While this approach has seen heuristic success, {it does not generally find a solution to \cref{eq:ell_0}. Many studies have addressed this gap by probabilistic guarantees \cite{candes2005decoding,chandrasekaran2012convex,vidyasagar2019introduction} such as the restricted isometry property (RIP).}

Despite extensions of (RIP) to structured random matrices (see, e.g., \cite{foucart2013invitation,rauhut2010compressive}), in several practical applications where $V$ is a structured matrix, these statistical guarantees may be inconclusive, intractable \cite{tillmann2013computational} or too conservative. {The latter is particularly evident in discrete-time optimal control problems, such as
\begin{equation}\label{eq:ctrl_intro}
	\begin{aligned}
		& \underset{u}{\text{minimize}}
		& & \sum_{i=0}^{N-1} |\sign(u(t))| \\
		& \text{subject to}
		& & x({t+1}) = Ax(t)+bu(t), \\ 
		& & & x(0) = \xi,\; x(N) = 0,\\
	\end{aligned}
\end{equation}
where $A\in\mathbb{R}^{m\times m}$, $b, \ \xi \in\mathbb{R}^{m \times 1}$ and $0 \leq t \leq N-1$. Zero values of $u$ in such problems can be interpreted as no (or least possible) fuel consumption. Since
\begin{equation*}
	x(N) = A^N x(0) + \begin{bmatrix} b & A b & \dots & A^{N-1} b\end{bmatrix} \begin{bmatrix}
		u({N-1})\\
		\vdots\\ 
		u(0)
	\end{bmatrix},
\end{equation*}
we can reformulate \cref{eq:ctrl_intro} into the form of \cref{eq:ell_0} by choosing 
\begin{equation}
	V = \begin{bmatrix} b & A b & \dots & A^{N-1} b\end{bmatrix}\quad \text{and} \quad y = -A^{N}x(0), \label{eq:recov_state_fuel}
\end{equation}
i.e., $V$ is an \emph{extended controllability matrix}. Even in cases where \cref{eq:ctrl_intro} has a unique solution, \cref{eq:ell_1} does not recover that solution (see \cref{fig:ctrl_results_intro} for an example).} \begin{figure}[t]
	\centering
	\begin{tikzpicture}
		\begin{axis}[xmin = 0, xmax = 39, xlabel={$t$}, grid, ylabel = {$u(t)$}, height = 5 cm, width = 12 cm]
			\addplot [ycomb, thin, mark=square*,mark options={fill=red,draw=black},mark size=1.5 pt]file{u_2.txt}; \label{line:recov_v2}

		\end{axis}
	\end{tikzpicture}
	
	\caption{{The solution of \cref{eq:ell_1} applied to \cref{eq:ctrl_intro} with $A = \diag(\begin{bmatrix} 0.8 & 0.7 & 0.6 & 0.5 & 0.4 \end{bmatrix})$, $b = \begin{bmatrix}
	    1 & \cdots & 1 
	\end{bmatrix}^\transp$, $\xi = \begin{bmatrix}
8.0632 & 33.9728 & 163.7151 & 1022 & 9534.2432
\end{bmatrix}^\transp$ and $N=40$. While the unique solution to \cref{eq:ctrl_intro} is given by $u^\ast$ with $u^\ast(0) = 1$, $u^\ast(9) = -1$ and $u^\ast(t) = 0$ for $t \not \in \{0,9\}$, \cref{eq:ell_1} fails to recover $u^\ast$ and instead returns \ref{line:recov_v2} with twice as many non-zero entries.}}
\label{fig:ctrl_results_intro}
\end{figure} Yet, for analogous continuous-time problems, the equivalence between corresponding $L_0$- and $L_1$-norm minimization problems is well established \cite{athans1966optimal,Nagahara2016handsoff}. One would, therefore, expect this equivalence to also hold when such problems are discretized and solved numerically, i.e., when considering \cref{eq:ell_0,eq:ell_1}. This observation indicates that control theory can offer new deterministic insights into the equivalence of such problems {by first investigating the affects of different $(A,b)$ and then generalizing them to similar instances, where $V$ is not necessarily associated with system matrices.} 

In this study, rather than focusing on conditions under which \cref{eq:ell_1} succeeds, {we exploit the above observation to develop a deterministic, control theory-inspired understanding of its failure cases. This may be used for a better design of sparse optimization problems, e.g., by the a priori selection of basis functions (features) or the strategic placement of poles, as well as to rule out whether a solution to \cref{eq:ell_1} fails to solve \cref{eq:ell_0}.} We, therefore, we derive a condition on $V$ under which a solution $u^\ast$ of \cref{eq:ell_0} cannot be recovered by \cref{eq:ell_1}. {It is shown that this occurs if there exists at least one non-zero entry $u^\ast_i \neq 0$ for which 
\begin{equation*}
\|(V_{(:,(1:m))})^{-1} V_{(:,\{i\})}\|_1 < 1,     
\end{equation*}
where $V_{(:,(1:m))}$ stands for the submatrix consisting of the first $m$ columns of $V$ and $V_{(:,\{i\})}$ is the $i$-th column of $V$.} Although identification of such column indices is generally intractable {(see, e.g., \cref{ex:hankel_multi_mod} in \Cref{sec:main})}, we leverage the framework of total positivity \cite{karlin1968total} to provide a class of non-trivial matrices $V$ for which this issue can be resolved. We show that if $V$ is the extended controllability matrix of a linear time-invariant discrete-time system and all its minors of order $m$ share the same strict sign, then the $\ell_1$-norm of its transformed columns forms a log-concave vector. Based on this property, we further show that if there exists a transformed column in $V$ for which our condition prohibits a non-zero entry in $u^\ast$, then the same holds for all subsequent columns in $V$ {(see~\cref{fig:unimodal_intro} for an example)}. 
\begin{figure}
	\centering
	\begin{tikzpicture}
		\begin{groupplot}[
			group style={
				group size=1 by 1,
				vertical sep=1.25 cm,
			},
			height=5 cm,
			width=12cm,
			]
			\nextgroupplot[
			xmin=1, xmax=50, 
			xlabel={$k$}, 
			ylabel={$p_k$},
			ytick={0, 20, 40, 60}, 
			ymax=60,
			grid,
			ytick={0.01,0.1,1,10,100},
			ymode=log
			]
			\addplot[line width=1.5pt, blue] 
			table [col sep=space,y expr=\thisrowno{1}]{unimodal_vector_v2.txt}; 
			\label{line:unimod_v2}
			
		\end{groupplot}
	\end{tikzpicture}
	
	\caption{{Illustration of $p_k := \normA{V_{(:,(1:m))}^{-1}V_{(:,\{k\})}}{\ell_1}$ for $V = \begin{bmatrix} b & A b & \dots & A^{49} b\end{bmatrix}$ with $A$ and $b$ as in \cref{fig:ctrl_results_intro}. Our results guarantee that $p_k$ is log-concave for $k > m$, which appears as a concave curve on the logarithmic y-axis, and therefore $p$ is unimodal (i.e., single-peaked). Consequently, since $p_k$ crosses $1$ at $k = 36$, $p_k$ must remain below $1$ for all $k \geq 36$. We further show that this implies that a solution $u^\ast$ to \cref{eq:ctrl_intro} can be recovered by \cref{eq:ell_1} only if $u^\ast(t) = 0$ for all $t \ge 36$, providing a tractable explanation for the failure of \cref{eq:ell_1} in \cref{fig:ctrl_results_intro}.}}
	\label{fig:unimodal_intro}
\end{figure}{Thus, our index search is mapped onto a simple bisection process of finding the smallest such index. In other words, any solution of \cref{eq:ell_1} that includes a non-zero entry beyond this critical index cannot be a solution to \cref{eq:ell_0}.}

{Our findings also generalize beyond extended controllability matrices by imposing minor assumptions onto the column-wise forward difference of general matrices $V$. These conditions are suitable for a compressed sensing setting, i.e., where $V$ and $y$ stem from sampling/measuring the response of a larger linear operator from which $V$ inherits these properties. It is important to note that our minor-based conditions can be verified in polynomial time (see, e.g., \cite{carter2021complexity}) and by that also provide a tractable way to guarantee that \cref{eq:ell_0} has a unique solution. The latter ensures that our failure guarantees are not limited to trivial non-unique instances under which any convex relaxation is prone to fail due to forming convex combinations.} 
	
	In case of the extended controllability matrix, we, further, present a condition on the poles of an underlying controllable system for which \cref{eq:ell_1} only manages to recover the solution \cref{eq:ell_0} in the trivial $m=n$ instance. 
Finally, our log-concavity result complements work on the log-concavity of multivariate symmetric polynomials such as \cite{sagan1992log}.  

{The remainder of the paper is organized as follows. \Cref{sec:setup} introduces the notation and theoretical basis for our analysis. In \Cref{sec:main}, we present our main results regarding the interconnection of sparse optimization and total positivity from an optimal control viewpoint. Building upon these theoretical insights, \Cref{sec:exm} provides further illustrative examples based on the well-studied fuel-optimal control problem \cref{eq:ctrl_intro}, and conclusions are drawn in \Cref{sec:cncl}. Proofs can be found in the appendix.}

\section{Preliminaries} \label{sec:setup}
\subsection{Notations}
\subsubsection{Sets}
{We use $\mathbb{R}_{\geq 0}$ to denote sets of non-negative reals (with analogous notations for strictly positive and reversed inequalities, as well as integers $\mathds{Z}$). The set of natural numbers is denoted by $\mathbb{N} = \mathds{Z}_{\geq 1}$ and for $k,l\in\mathbb{N}$ we define 
	$(k:l) := \{k,k+1,\ldots,l\}$ if $k\leq l$ and equal to the empty set otherwise. Finally, the collection of strictly increasing \(r\)-tuples with values in $(1:n)$ is denoted by} $$
\mathcal{I}_{n,r}
:=\bigl\{v=(v_1,\dots,v_r)\in\mathbb{N} : 1\le v_1<\cdots<v_r\le n\bigr\}.$$

\subsubsection{Matrices}
For a matrix $X = (x_{ij}) \in\Rmn$, we use $\text{Im}(X)$ to denote its image and say that $X$ is \emph{non-negative} ($X \geq 0$) if all $x_{ij} \geq 0$, and \emph{positive} ($X > 0$) if all $x_{ij} > 0$. The \emph{Moore–Penrose inverse} of $X$ is denoted by $X^\dagger$. 
Submatrices of $X$ that are formed from row indices $\mathcal{I} \subset (1:m)$ and column indices $\mathcal{J} \subset (1:n)$ are denoted by $X_{(\mathcal{I},\mathcal{J})}.$ For brevity, we use the colon $:$ notation to denote $X_{(:,\mathcal{J})} = X_{((1:m),\mathcal{J})}$ and $X_{(\mathcal{I},:)} = X_{(\mathcal{I},(1:n))}$. These notions are also used for $\mathcal{I} \in \mathcal{I}_{m,r}$ and $\mathcal{J} \in \mathcal{J}_{n,p}$. A \emph{(consecutive) $j$-minor} of $X$ is a minor which is constructed of $j$ columns and $j$ rows of $X$ (with consecutive indices). We also use these notations for matrices with countably infinite rows and/or columns. If $X \in \mathds{R}^{n \times n}$, we denote its spectrum by $\sigma(X) = \{\lambda_1(X),\dots,\lambda_n(X)\}$ and by $I_n$ the identity matrix in $\mathds{R}^{n \times n}$. Diagonal $X$ is abbreviated by $X = \diag(v)$ with $v \in \mathds{R}^{n}$ satisfying $ v_i = x_{ii}$. Further, $\mathbf{1}_n \in \mathds{R}^n$ is the all ones vector and $e_{i} \in \mathds{R}^n$ the $i$-th canonical unit vector. We, further, make use of what we call the \emph{row forward difference} $\Delta X \in \mathds{R}^{n-1 \times n}$ of $X \in\mathbb{R}^{n\times m}$, where $\Delta = (\delta_{ij})\in \mathds{R}^{n -1 \times n}$ is defined by
\begin{equation*}
	(\delta)_{ij} = \begin{cases}
		-1 & \text{if } j = i\\
		1& \text{if }j = i+1 \\
		0 & \text{else}
	\end{cases}.
\end{equation*}
Further, the anti-diagonal matrix $K_m = ((k_m)_{ij}) \in\mathbb{R}^{m\times m}$ is defined by 
\begin{equation*}
	(k_m)_{ij} = \begin{cases}
		(-1)^{j-1} & \text{if } i+j = m+1\\
		0 & \text{else}
	\end{cases}.
\end{equation*}
\subsubsection{Functions}
The \emph{(0-1)-indicator function} of a subset $\mathbb{A} \subset \mathds{R}$ is defined by
\begin{equation*}
	\mathbb{1}_{\mathbb{A}}(x) := 
	\begin{cases}
		1, & \text{if } x \in \mathbb{A},\\[1mm]
		0, & \text{otherwise.}
	\end{cases}
\end{equation*}
and the \emph{floor function} by $\lfloor \cdot \rfloor$. The set of \emph{$k$-time continuously differentiable scalar functions} on the open interval $(a,b)$ are denoted by $C^{k}(a,b)$. 

\subsection{Dual Minimum Norm Problem}
Given a norm $\normA{\cdot}{}$ on $\mathds{R}^n$, its \emph{dual norm} is defined by $\normA{x}{\ast} := \sup_{\normA{d}{}\leq 1} x^\transp d$. A standard example (see \cite{luenberger1968optimization}) is $$\normA{x}{\ell_1 \ast} = \normA{x}{\ell_\infty} :=  \max_{i \in (1:n)} |x_{i}|.$$ A vector $u\in \mathbb{R}^n$ is said to be \emph{aligned} with $y\in  \mathbb{R}^n$ if $u^\transp y=\|y\|_\ast\|u\|$. %
The following so-called \emph{dual approximation theorem} \cite[Corollary~5.9.1]{luenberger1968optimization}, will be essential in our analysis. 

\begin{prop}\label[proposition]{lem:opt_cond}
	Let $V \in\mathbb{R}^{m\times n}$ and $y\in\mathbb{R}^m$ be such that $$D: =\{u\in \mathbb{R}^n:Vu=y\} \neq \emptyset.$$ Then,
	$\min_{u\in D}\normA{u}{} = \max_{\|V^\transp\beta\|_\ast\leq 1} y^\transp \beta$ and if the minimum is achieved by $u^\ast$ and the maximum by $\beta^\ast$, then $u^\ast$ is aligned with $V^\transp \beta^\ast$.
\end{prop}

\subsection{Linear Time-Invariant Discrete-Time Systems}
{We consider \emph{causal linear time-invariant (LTI) discrete-time system} (see, e.g., \cite{antsaklis2006linear,rugh1996linear,kwakernaak1991modern})} defined by $A \in \mathds{R}^{m \times m}$ and $b, c^\transp \in \mathds{R}^m$ as
\begin{equation}\label{eq:LTI}
	\begin{aligned}
		x(t+1) &= A x(t) + b u(t), \\
		y(t) &= c x(t).
	\end{aligned}
\end{equation}
where we refer to $u(t), y(t)\in \mathds{R}$, and $x(t) \in \mathds{R}^m$ as the input, output, and state, respectively. Given initial conditions $x(0)$, the solution at $t=N \geq 0$ is given by 
\begin{equation}\label{eq:sol_state}
	x(N) = A^N x(0) + \mathcal{C}^N(A,b) \begin{bmatrix}
		u({N-1})\\
		\vdots\\ 
		u(0)
	\end{bmatrix}
\end{equation}
where $$\mathcal{C}^N(A, b) = \begin{bmatrix} b & Ab & A^2 b & \dots & A^{N-1}b \end{bmatrix}$$
is the so-called \emph{extended controllability matrix}. The \emph{extended observability matrix} is denoted by $\mathcal{O}^{N}(A,c):= \mathcal{C}^{N}(A^\transp,c^\transp)^\transp$
and in case of $N = \infty$, we simply write $\mathcal{C}(A,b)$ and $\mathcal{O}(A,c)$, respectively. The pair $(A,b)$ is called \emph{controllable}, if $\det(\mathcal{C}^{m}(A,b)) \neq 0$ and \emph{observable}, if $\det(\mathcal{O}^{m}(A,c)) \neq 0$. 

The \emph{impulse response} of \cref{eq:sol_state} is given by $g(t) = cA^{t-1}b \mathds{1}_{\mathds{Z}_{>0}}$ and its \emph{transfer function} by 
\begin{equation*}
	G(z) = \sum_{t=1}^{\infty} g(t)z^{-t}.
\end{equation*}
The triplet $(A,b,c)$ is also called a \emph{realization} of $G(z)$, which is \emph{minimal}, if $(A,b)$ and $(A,c)$ are controllable and observable, respectively. Finally, let us define the Hankel matrices corresponding to \cref{eq:sol_state} by
\begin{equation*}
	H_{M,N}(A,b,c) := \mathcal{O}^M(A,c) \mathcal{C}^N(A,b) = \begin{bmatrix}
		g(1) & \dots & g(N)\\
		\vdots & & \vdots\\
		g(M) & \dots & g(N+M-1)
	\end{bmatrix}.  
\end{equation*}

\subsection{Variation bounding and total positivity}
\label{sec:vardim}
Variation bounding provides key insights into how sign changes evolve under linear transformations \cite{karlin1968total}. Next, we introduce formal definitions and relevant properties essential for our analysis.
\begin{defn}
	The \emph{variation} of a vector $u\in\mathds{R}^n$ is defined as
	\begin{equation*}
		\vari{u} := \sum_{i = 1}^{p-1} \mathbb{1}_{\mathbb{R}_{< 0}}(\tilde{u}_i \tilde{u}_{i+1}) \quad \text{and}  \quad 
		\vari{0} := -1
	\end{equation*}
	where $\tilde{u} \in \mathds{R}^s$, $s \leq n$, is obtained by removing all zero entries from $u$.
\end{defn}
\begin{defn}\label{def:vb_k}
	$X \in\mathbb{R}^{n\times m}$ is said to be \emph{$k$-variation bounding} ($X \in \vb{k}$), $0 \leq k \leq \min\{m-1,n-1\}$, if for all $u \in \mathds{R}^m$ with $\vari{u} \leq k$, it holds that $\vari{Xu} \leq k$.
\end{defn}
To characterize $X \in \vb{k}$, we require the so-called \emph{$r$-th multiplicative compound matrix} {(see, e.g., \cite{wedderburn1934lectures,karlin1968total,horn2012matrix,aitken1939determinants,fallat2011totally,barshalom2022compound,pinkus2009totally,cullis1925matrices} for a historical account)}: for $X \in \mathbb{R}^{n \times m}$ and $r \in (1:\min \{m,n\})$, the $(i,j)$-th entry of the $r$-th multiplicative compound matrix $X_{[r]} \in \mathbb{R}^{\binom{n}{r} \times \binom{m}{k}}$ is defined via $\det(X_{(I,J)})$, where $I$ and $J$ are the $i$-th and $j$-th element of the $r$-tuples in $\mathcal{I}_{n,r}$ and $\mathcal{I}_{m,r}$, respectively, assuming \emph{lexicographical ordering}. For example, if $X \in \mathbb{R}^{3 \times 3}$, then
        \[ X_{[2]}=\begin{bmatrix}
	\det(X_{(\{1,2\},\{1,2\})}) & \det(X_{(\{1,2\},\{1,3\})}) & \det(X_{(\{1,2\},\{2,3\})}) \\
	\det(X_{(\{1,3\},\{1,2\})}) & \det(X_{(\{1,3\},\{1,3\})}) & \det(X_{(\{1,3\},\{2,3\})}) \\
	\det(X_{(\{2,3\},\{1,2\})}) & \det(X_{(\{2,3\},\{1,3\})}) & \det(X_{(\{2,3\},\{2,3\})})
\end{bmatrix}.
\]
The following properties of the multiplicative compound matrix {(see, e.g., \cite[Subsections~0.8.1 and~0.4.4]{horn2012matrix}, \cite[Chapter~V, Theorem 1]{wedderburn1934lectures} or \cite[Section~38]{aitken1939determinants})} will be  used:

\begin{lem}\label[lemma]{lem:comp_properies}
	Let $X \in \mathbb{R}^{m \times n}$ and $Y \in \mathbb{R}^{n \times p}$. Then, 
	\begin{enumerate}[label=\roman*.]
		\item {$X_{[r]}Y_{[r]} = (XY)_{[r]}$ for any $r \in (1:\min\{m,n,p\})$}.
		\item $\rk(X) = \min \{m,n\}$ if and only if $X_{[\min \{m,n\}]} \neq 0$.
	\end{enumerate}
\end{lem}
Under a rank assumption, $X \in \vb{k}$ can be characterized by the notions of $k$-sign consistency and total positivity (see~\cite{roth2024system} for a recent discussion). 
{\begin{defn}
		Let $X\in\mathbb{R}^{m\times n}$ and $k\leq\min\{m,n\}$. Then,
		\begin{itemize}
			\item $X$ is called (strictly) $k$-sign consistent ($X \in\text{(S)SC}_k$) if $X_{[k]} \geq 0$ or $X_{[k]} \leq 0$ (respectively, $X_{[k]} > 0$ or $X_{[k]} < 0$). 
			\item $X$ is called (strictly) $k$-totally positive  ($X\in\text{(S)TP}_k$) if $X_{[j]} \geq 0$ (respectively, $X_{[j]} > 0$) for all $j\in (1:k)$. If $k=\min\{m,n\}$, then $X$ is called (strictly) totally positive.
		\end{itemize}
		These notations are also used for the operators $\mathcal{C}(A,b)$ and $\mathcal{O}(A,c)$, meaning that $\mathcal{C}^N(A,b)$ and $\mathcal{O}^N(A,c)$, respectively, fulfill these properties for all $N \geq k$. 
	\end{defn}
	\begin{lem}\label[lemma]{lem:vb_m_sc_m}
		Let $X \in \mathds{R}^{n \times m}$ and $k < \rk(X)$ be such that any $k$ columns of $X$ are linearly independent. Then, $X \in \vb{k-1}$ if and only if $X \in \text{SC}_k$. Moreover, $X \in \text{VB}_{m-1}$, $n > m = \rk(X)$, if and only if $X \in \text{SC}_m$.  
	\end{lem}
}
It is important to note that sign consistency and total positivity can often be verified efficiently in polynomial time, e.g., by the following propositions (see~\cite[Theorem~2.3.1]{karlin1968total} and \cite[Proposition~5]{grussler2020variation}):
\begin{prop}\label[proposition]{prop:consecutive_karlin}
	Let $X \in \Rmn$ and $n \geq m$ be such that
	\begin{enumerate}[label=\roman*.]
		\item all consecutive $s$-minors of $X_{(:,1:s)}$ share the same strict sign for $s \in (1:m-1)$
		\item all consecutive $m$-minors of $X$ share the same (strict) sign.
	\end{enumerate}
	Then, $X \in \text{(S)SC}_m$.
\end{prop}
\begin{prop}
	\label[proposition]{prop:consecutive_old}
	Let $X \in \Rmn$ and $k \leq \min\{m,n\}$ be such that
	\begin{enumerate}[label=\roman*.]
		\item all consecutive $j$-minors of $X$ are 
		strictly positive, $j \in (1:k-1)$,
		\item all consecutive $j$-minors of $X$ are
		nonnegative (strictly positive).
	\end{enumerate}
	Then, $X \in \text{(S)TP}_k$. In the strict case, these are also necessary.  
\end{prop}
{By application of repeated Dodgson condensation as in \cite{carter2021complexity}, it can be shown that the total complexity of verifying the conditions of these propositions is bounded by $\mathcal{O}(m^2n)$ if $m \leq n$. A geometric interpretation and generalization of these propositions has recently been discussed in \cite{grussler2026discrete}.} \cref{prop:consecutive_karlin,prop:consecutive_old} may simplify further in case of $\mathcal{C}^N(A,b) \in \text{(S)TP}$ \cite[Theorem~3.5]{grussler2021internally}. Using these results, it is easy to check that if $A = \diag(x)$, $x \in \mathds{R}^m$, and $b = \mathbf{1}_m$, then the so-called \emph{Vandermonde matrix}
\begin{equation*}
	\mathcal{C}^N(A,b) =  \begin{bmatrix}
		1 & x_1 & \dots & x_1^{N-1}\\
		\vdots & \vdots & & \vdots \\
		1 & x_m & \dots & x_m^{N-1}
	\end{bmatrix} \in \text{(S)TP}_m
\end{equation*}
if and only if $0 (<) \leq x_1 \leq x_2 \leq \dots \leq x_m$. (Strict) total positivity can also be used to parametrize (strictly) sign-consistent matrices as follows \cite[Proposition~2.1]{pena}:
\begin{lem}\label[lemma]{lem:pena}
	Let $X \in \mathds{R}^{n \times m}$ , $n > m$, be such that $\det (X_{(1:m,:)}) \neq 0$ and
	\begin{equation*}
		X(X_{((1:m),:)})^{-1}K_m = \begin{bmatrix}
			K_m \\
			C
		\end{bmatrix}.
	\end{equation*}
	Then, there is a bijection between the $m$-minors $\det(X_{(\gamma,:)})$ for all $\gamma \in \mathcal{I}_{n,m} \setminus (1:m)$ and the entries of $C_{[m]}$. In this bijection, the corresponding determinants coincide up to the sign of $\det(X_{(1:m,:)})$.
\end{lem}
In other words, $X \in \mathds{R}^{n \times m}$, $n > m$, is (strictly) $m$-sign consistent if and only if the matrix $C$ in \cref{lem:pena} is (strictly) totally positive. 
{Finally, a vector $a \in \mathbb{R}^n$ with $a \geq 0$ is said to be \emph{log-concave} if it satisfies 
	\begin{equation*}
		a_{k+1}^2 \geq a_{k}a_{k+2} \text{ for all }k \geq 2,
	\end{equation*}
	and \emph{unimodal} if there exists an $l \in (1:n)$ such that 
	\begin{equation*}
		a_1 \leq \ldots \leq a_l \geq \ldots \geq a_n.	
	\end{equation*}
	We also use these notions in the infinite case of sequences $\{a_k\}_{k \geq 1}$. Equivalently, $a$ is unimodal if and only if $\vari{\Delta a} \leq 1$ and the possible sign switch occurs from positive to negative. In particular, log-concavity is known to imply unimodality \cite{karlin1968total}.}

\section{Main results} \label{sec:main}
Our work establishes tractable conditions under which the basis pursuit approach fails to solve \cref{eq:ell_0}. We first derive a general condition based on the $\ell_1$-norm of the columns of a general matrix $V$, and demonstrate its relevance to structured matrices. A uniqueness guarantee for the solution of \cref{eq:ell_0} shows that our conditions are not limited to trivial non-unique cases, where basis pursuit inherently fails. We then show that both, the failure and uniqueness conditions can be efficiently verified for system matrices such as $V = \mathcal{C}^N(A,b) \in \text{SC}_m$. Leveraging total positivity theory \cite{karlin1968total,roth2024system}, we further extend our findings to general $V \in \text{SC}_m$ and discuss its suitability for compressed sensing and connection to geometric modeling.

\subsection{General Failure Condition}
\begin{thm}\label{thm:l1_necessary}
	Let $u^\ast \in \mathds{R}^n$ be a solution to \cref{eq:ell_0}, where $V \in \mathds{R}^{m \times n}$, $n \geq m$, is such that
	\begin{enumerate}[label=\roman*.]
		\item $\rk(V_{(:,(1:r))}) = \rk(V) = r$.
		\item  $V = V_{(:,(1:r))} V_{(:,(1:r))}^\dagger V$.
	\end{enumerate}
	 Then, $u^\ast$ is not a solution of \cref{eq:ell_1} if there exists an $i^\ast \in (r+1:n)$ with $u^\ast_{i^\ast}\neq 0$ and
	\begin{equation}
		\|V_{(:,(1:r))}^\dagger V_{(:,\{i^\ast\})}\|_{\ell_1} < 1. \label{eq:colum_cond}
	\end{equation}
\end{thm}
Note that our condition on $V$ merely requires that all its columns lie within the image of $V_{(:,(1:r))}$. This condition is trivially satisfied if $\det(V_{(:,1:m)}) \neq 0$. It also arises naturally, e.g., if $V = \mathcal{C}^N(A,b)$ or if we can factor $V = L R$ with $L \in \mathds{R}^{m\times r}$, $R \in \mathds{R}^{r \times n}$ such that $\det(R_{(:,1:r)}) \neq 0$ as in the case of Hankel matrices 
$$H_{m,n}(A,b,c) = \mathcal{O}^{m}(A,c)\mathcal{C}^{n}(A,b),$$
where the realization $(A,b,c)$, $A \in \mathds{R}^{r \times r}$, can be assumed to be a minimal and $m, n \geq r$. 
The following is a simple illustration of \cref{thm:l1_necessary}. 
\begin{exm}\label[example]{ex:simple}
	Let $V=\mathcal{C}^4(A,b)\in \mathbb{R}^{3\times 4}$ with $A=\diag(\begin{bmatrix}
		0.17 & 0.23 & 0.4
	\end{bmatrix})$ and $b=\mathbb{1}_3$. Then, since
	\begin{equation*}
		\|V_{(:,(1:3))}^{-1} V_{(:,\{4\})}\|_{\ell_1} = \|\begin{bmatrix}
			0.0143&-0.1877&0.78
		\end{bmatrix}\|_{\ell_1} < 1,
	\end{equation*}
	it follows from \cref{thm:l1_necessary} that a solution $u^\ast$ to \cref{eq:ell_0} cannot be recovered by solving \cref{eq:ell_1} if $u^\ast_4 \neq 0$.
\end{exm}
{
	\begin{rem}
		\cref{thm:l1_necessary} can be utilized as follows: columns with index $i^\ast$ fulfilling \cref{eq:colum_cond} can be removed before solving \cref{eq:ell_1} as they are irrelevant to our attempt of solving \cref{eq:ell_0} via \cref{eq:ell_1}. A goal of this work is to investigate matrix structures for which such indices can be computed efficiently.  
\end{rem}}
Next, let us note that our failure guarantees do not merely correspond to the trivial instances where \cref{eq:ell_0} admits more than one solution, i.e., \cref{eq:ell_1} fails due to the formation of convex combinations. {To this end, we use the following uniqueness guarantee \cite{gorodnitsky1997sparse}}:
\begin{prop}\label[proposition]{thm:uniquness}
	Let $V \in\mathbb{R}^{m\times n}$ and $y \in \mathds{R}^m$ be such that:
	\begin{enumerate}[label=\roman*.]
		\item Any $m$ columns of $V$ are linearly independent. \label{item:lin_ind_m}
		\item $\exists u^\ast \in \mathds{R}^n:\; V u^\ast = y$ with $\normA{u^\ast}{\ell_0} \leq \lfloor m/2\rfloor$. \label{item:sparse_bound}
	\end{enumerate}
	Then, $u^\ast$ is the unique solution to \cref{eq:ell_0}. In particular, if $u^\ast$ is a solution to \cref{eq:ell_1}, it also solves \cref{eq:ell_0}. 
\end{prop}
In \cref{ex:simple}, it is easy to verify, using \cref{lem:comp_properies} and the fact that $V\in\text{STP}$ is a Vandermonde matrix, that any $m$ columns of $V$ are linearly independent. Thus, even when \cref{thm:uniquness} guarantees that \cref{eq:ell_0} has a unique solution, basis pursuit may fail to recover it. 
\begin{rem}
	\cref{lem:pena,prop:consecutive_old} imply that \cref{item:lin_ind_m} in \cref{thm:uniquness} can be verified numerically efficiently, i.e., in polynomial time, if $V \in \text{SC}_m$ \cite{roth2024system}. A goal of this work is to show that the $\text{SC}_m$ property also facilitates a tractable verification of \cref{thm:l1_necessary}. 
\end{rem}
Next, let us define 
\begin{equation}
	p := (\normA{V_{(:,(1:r))}^{\dagger }V_{(:,\{i\})}}{\ell_1})_{i \in (1:n)}\in \mathds{R}^n \label{eq:p_unimod}
\end{equation}
where the condition in \cref{thm:l1_necessary} corresponds to entries $p_{i^\ast} < 1$, $i^\ast \in (r+1:n)$ and $p_{(1:r)} = \mathbf{1}_r$, where $r = \rk(V_{(:,(1:r))}) = \rk(V)$. In general, $p$ may contain many local extrema as demonstrated by the following example.
\begin{exm}\label[example]{ex:hankel_multi_mod}
	Let $V = H_{5,N}(A,b,c)$ for $N \geq 20$ and any realization $(A,b,c)$ of the transfer function $$G(z) = \frac{1}{(z-0.8)(z-0.5)(z-0.1)}.$$ Then, $p$ in \cref{eq:p_unimod} drops below one at $k=11$ and again for $k \ge 17$ (see~\cref{fig:p_multi_mod}). While by the asymptotic stability of $G(z)$, $p_k$ has to eventually drop below $1$ and remain there for sufficiently large $k$ and $N$, the transient behavior of $p$ makes an efficient determination of all entries $p_k < 1$ challenging. 
\end{exm}
\begin{figure}[t]
	\centering
	\begin{tikzpicture}
		\begin{groupplot}[
			group style={
				group size=1 by 2,
				vertical sep=1.25 cm,
			},
			height=5 cm,
			width=12cm,
			]
			
			\nextgroupplot[xmin=1, xmax=20, 
			xlabel={$k$}, 
			ylabel={$p_k$},
			grid,
			ymode = log,
			ytick={0.01, 0.1, 1,10,100}, 
			]
			\addplot[line width=1.5pt, blue] 
			table [col sep=space,y expr=\thisrowno{1}]{p_multi_mode.txt}; 
			
		\end{groupplot}
	\end{tikzpicture}
	
	\caption{{Illustration of the vector $p := (\normA{V_{(:,1:m)}^{\dagger }V_{(:,\{k\})}}{\ell_1})_{k \in (1:n)}\in \mathds{R}^N$ corresponding to \cref{ex:hankel_multi_mod} with $N = 20$. 
		The vector $p$ exhibits many local extrema and drops below $1$ at $k = 11$ and for $k \in (17:20)$, which makes the application of a principled search method challenging.}
	}
	\label{fig:p_multi_mod}
\end{figure}
\subsection{Tractability for Systems Matrices}
In this subsection, we want to show that by restricting ourselves to system matrices/operators, it is possible to drive tractable conditions for determining the indices with $p_k < 1$. The first class of matrices we examine is \(V=\mathcal{C}^N(A,b)\), which satisfies \(p_{m+i}<1\) for all \(i\ge 1\). 
To keep our derivations independent of the number of columns in $V$, we present our result in the infinite-dimensional case. 
\begin{thm}\label{thm:char_poly}
	Let $A\in\mathbb{R}^{m \times m}$, $b \in \mathds{R}^m$ and $V = \mathcal{C}(A,b)$ and $r =\rk(V)$. If the characteristic polynomial 
	\begin{equation*}
		\det(sI-V_{(:,(1:r))}^\dagger A V_{(:,(1:r))}) = s^r + \alpha_{r-1} s^{r-1} + \dots + \alpha_1 s + \alpha_0 
	\end{equation*}
fulfills $\sum_{i=0}^{r-1} |\alpha_i| < 1$, then 
\begin{equation*}
	\forall i^\ast > r: \; \normA{V_{(:,(1:r))}^{\dagger }V_{(:,\{i^\ast\})}}{\ell_1} < 1.
\end{equation*}
Consequently, by \cref{thm:l1_necessary}, if $V=\mathcal{C}^N(A,b)$ in \cref{eq:ell_1}, then $u^\ast \in \mathds{R}^N$ is not a solution of \cref{eq:ell_1} if there exists an $i^\ast \in (r+1:N)$ with $u^\ast_{i^\ast}\neq 0$. 
\end{thm}
\cref{thm:char_poly} has the advantage that it can be verified efficiently and is independent of the number of columns in $\mathcal{C}^N(A,b)$ as well as almost all $b$. In particular, it shows that if $(A,b)$ is controllable then \cref{eq:ell_1} succeeds if and only if $N = m$, i.e., when no optimization is required. This drastic failure may in practical applications be prevented by a better design of $A$. Note that in case of $\rk(V) < m$, $(V_{(:,(1:r))}^\dagger A V_{(:,(1:r))},V_{(:,(1:r))}^\dagger b)$ simply corresponds to the dynamics after projecting onto the controllable subspace, i.e., the uncontrollable modes are irrelevant to our ability to find, e.g., a sparest control signal as discussed in \Cref{sec:exm}. As the following simple example shows, one can cannot conclude that if the characteristic polynomial of $A$ satisfies the coefficient constraint in \cref{thm:char_poly} that the same is true after the projection:
\begin{exm}
	Let $V = \mathcal{C}(A,b)$ with 
	\begin{equation*}
		A = \begin{bmatrix}
			0 & -0.7 & 0\\
			1 & -0.5 & 0\\
			0 & 0 & 0.5
		\end{bmatrix}\; \quad b = \begin{bmatrix}
		1\\
		1\\
		0
		\end{bmatrix}
	\end{equation*}
	Then, the characteristic polynomial of $A$
	\begin{equation*}
		\det(sI-A) = s^3 + 0.45 s - 0.35
	\end{equation*}
	fulfills coefficient constraint in \cref{thm:char_poly}. However, as $\rk(V) = 2$, we need to consider
	\begin{equation*}
			\det(sI-V_{(:,(1:2))}^\dagger A V_{(:,(1:2))}) = s^2 +0.5s + 0.7,
	\end{equation*}
	which violates this constraint. 
\end{exm}
Unfortunately, \cref{thm:char_poly} is inconclusive for many pairs $(A,b)$ such as the one in \cref{ex:simple}. In order to mitigate this issue as well as the one observed in \cref{ex:hankel_multi_mod}, we want to derive conditions for when $p$ is unimodal. This property is desirable, since if $p_{i^\ast} < 1$ then $p_{i} < 1$ for all $i > i^\ast$, thereby reducing the verification of \cref{thm:l1_necessary} to a simple binary search problem. Since $p_{(1:r)} = \mathbf{1}_r$, it suffices to show unimodality of the vector $(\normA{V_{(:,1:m)}^{-1}V_{(:,\{i\})}}{\ell_1})_{i \in (r+1:n)}.$
We start with the case of $V = \mathcal{C}(A,b) \in \text{SC}_m$, where in order to analyze it, we first need to establish the following two lemmas:
\begin{lem}
	\label[lemma]{lem:spectrum}
	Let $A\in\mathbb{R}^{m \times m}$ and $b \in \mathds{R}^m$ be such that $(A,b)$ is controllable and $V= \mathcal{C}(A,b) \in \text{SC}_m$. Then $\sigma(A) \subset \mathds{R}_{\geq 0}$. 
\end{lem}  
\begin{lem}\label[lemma]{lem:pseudo_scm}
Let $X \in \mathds{R}^{n \times m}$, $n \geq m$, be such that $\rk(X) = m$. Then, $X \in \text{(S)SC}_m$ if and only if $X^\dagger \in \text{(S)SC}_m$ and the sign of the $m$-minors in $X$ and $X^\dagger$ coincide. 
\end{lem}
Based on these two lemmas, we can now arrive at our next main result:
\begin{thm}\label{thm:log_conv}
	Let $A\in\mathbb{R}^{m \times m}$, $b \in \mathds{R}^m$ and $V= \mathcal{C}(A,b)$ be such that $\rk(V) = r$ and $V \in \text{SC}_r$. Further, let 
\begin{equation*}
	W_k := K_r^\transp V_{(:,(1:r))}^{\dagger} V_{(:,\{r+k\})}, \; k \geq 1.
\end{equation*}
	Then, the sequences $\{g_{(i)}(k)\}_{k\ge1}$, $i \in (1:r)$, and $\{p(k)\}_{k\ge1}$ defined by 
	\begin{equation*}
		g_{(i)}(k) := W_k^\transp e_i \quad \text{and} \quad p(k) := \sum_{i=1}^r g_{(i)}(k) = W_k^\transp \mathbf{1}_r = (\normA{V_{(:,(1:r))}^{\dagger}V_{(:,\{r+k\})}}{\ell_1})
	\end{equation*}
	are log-concave. 
\end{thm}
As log‐concavity implies unimodality, \cref{thm:log_conv} provides conditions under which $V = \mathcal{C}^N(A,b)$ yields the desired unimodal $p$ defined in \cref{eq:p_unimod}. Note by \cref{lem:spectrum}, if $(A,b)$ is controllable, then $\mathcal{C}(A,b) \in \text{SC}_m$ requires that $\sigma(A) \subset \mathds{R}_{\geq 0}$. Otherwise, the result states that the uncontrollable modes are irrelevant, which is natural as these should not influence our ability to find, e.g., a sparsest control signal as discussed in \Cref{sec:exm}. Moreover, since $p$ is the impulse response of a system associated with $(A,b)$, $p$ will then necessarily contain entries strictly smaller than $1$ if $\sigma(A) \subset [0,1)$ and $N$ is sufficiently large. 

\begin{rem}
	In the special case $A = \diag(x)$ with $0 \leq x_1 \leq \cdots \leq x_m$ and $b = \mathbf{1}_m$, $V$ becomes a totally positive Vandermonde matrix. Here the sequences $g_{(i)}$ and $p$ consist of symmetric polynomials in $x_1,\dots,x_m$. For example, $$ g_{(1)}(k) = \sum_{1 \leq i_1 \leq i_2 \leq \dots \leq i_k \leq m} x_{i_1} x_{i_2} \cdots x_{i_k},$$ which for $k \in (1:m)$ is the so-called \emph{complete homogeneous symmetric polynomial of degree $k$ in $m$ variables} \cite{sagan1992log}. Our results are, therefore, also of independent interest, as they establish the log-concavity of such polynomial sequences for any fixed nonnegative vector $x \in\mathds{R}^m$. While these findings complement related work \cite{sagan1992log}, to the best of our knowledge no prior work has presented these results.
\end{rem}
Next, since impulse responses are invariant under similarity transformation, the same can be expected for \cref{thm:log_conv}.

\begin{cor}\label[corollary]{cor:transform}
	Let $A \in \mathds{R}^{m \times m}$ and $b \in \mathds{R}^m$ be as in \cref{thm:log_conv}. Further, let $T \in \mathds{R}^{m \times m}$ be invertible, ${V} = T \mathcal{C}(A,b)$ and define
	\begin{align*}
		g_{(i)}(k) &:=  |e_i^\transp V_{(:,(1:r))}^{\dagger}V_{(:,\{r+k\})}|, \;  i \in (1:r), \\
		p(k) &:= |\mathbf{1}_r^\transp V_{(:,(1:r))}^{\dagger}V_{(:,\{r+k\})}|, \; k \geq 1.
	\end{align*} 
	Then, $\{g_{(i)}(k)\}_{k\ge1}$ for each $i\in (1:r)$ and $\{p(k)\}_{k\ge1}$ are log-concave sequences.  
\end{cor}
For example, \cref{cor:transform} applies to Hankel matrices $$H_{m,n}(A,b,c) = \mathcal{O}^{m}(A,c)\mathcal{C}^{n}(A,b), \; n > m$$
with observable $(A,c)$ or the related Page matrices \cite{page}, which have recently gained attention in data‐driven control approaches (see, e.g., \cite{coulson2019regularized}), where \(\ell_1\)-norm regularization is employed to enhance robustness.

\subsection{General Tractability \& Compressed Sensing}
While \cref{thm:log_conv,cor:transform} give us conditions for retrieving our desired unimodality property, its structure is not very generic. For example, we would not be able to apply these results to the Hankel matrix cases $$V= \mathcal{O}^k(A,c) \mathcal{C}^{n}(A,b), \; k < m$$ or when $V$ is not related to an extended controllability matrix at all. Fortunately, as we only want $p$ to be unimodal --- not necessarily log-concave --- we can generalize our results via the notion of variation bounding. 
\begin{thm}\label{thm:unimodal}
	Let $V\in\mathbb{R}^{m\times n}$, $m < n$, be such that $V^\transp \in \text{SC}_{m}$, $\Delta (V^\transp) \in \text{VB}_{m-1}$ and $\det(V_{(:,1:m)}) \neq 0$. Then, $p \in \mathds{R}^{n}$ defined by $
	p_k := \normA{V^{-1}_{(:,(1:m))} V_{(:,\{k\})}}{
		\ell_1
	}$, $k \in (1:n)$, is unimodal.
\end{thm}

\begin{rem}
	By \cref{lem:vb_m_sc_m}, the requirement that $V^\transp \in \text{SC}_m$ can also be translated into $V^\transp \in \text{VB}_{m-1}$. Similarly, if $V^\transp$ does not contain any column that is a multiple of $\mathbf{1}_n$, then $\rk(\Delta (V^\transp)) = m$ and we can equivalently demand that $\Delta (V^\transp) \in \text{SC}_m$. In case that $\rk(V) = r$ and $V = V_{(:,(1:r))}V_{(:,(1:r))}^\dagger V$, one can apply \cref{thm:unimodal} to the projection $V_{(:,(1:r))}^\dagger V$. 
\end{rem}	
In the following, we present an illustrative example to demonstrate the benefit of \cref{thm:unimodal}.
\begin{exm}
	Let $A = \diag(\begin{bmatrix}
		0.1 & 0.2 & 0.3 & 0.4
	\end{bmatrix})$, $b = c^\transp =  \mathbf{1}_4$ and $$
	V = H_{3,6}(A,b,c) = \mathcal{O}^3(A,c)\mathcal{C}^6(A,b) \in \mathds{R}^{3 \times 6}.$$
	Since $\mathcal{O}^3(A,c)$, $\mathcal{C}^6(A,b) \in \text{STP}_3$, it follows from \cref{lem:comp_properies} that also $V \in \text{STP}_3$. Moreover, it is possible to directly verify that $\Delta (V^\transp) \in \text{SC}_3$ and $\rk(\Delta (V^\transp)) = 3$. Hence, by \cref{thm:unimodal}, we can conclude that $p$ defined in \cref{eq:p_unimod} is unimodal, despite $\mathcal{O}^3(A,c)$ being a non-square matrix. 
\end{exm}
Next we want to examine how the conditions in \cref{thm:unimodal} are suitable in a compressed sensing setting (see, e.g.,  \cite{rauhut2010compressive,vidyasagar2019introduction,brunton2022data}), where $V$ and $y$ in \cref{eq:ell_0} stem from sampling/measuring the response of a larger, possibly even infinite-dimensional, linear operator $u \mapsto Xu$. For example, $X$ may represent a dictionary of basis functions $f_1,\dots, f_n: \mathcal{T} \subset \mathds{R} \to \mathds{R}$ such as polynomials, wavelets, etc. (see, e.g., \cite{rauhut2010compressive,brunton2022data}), in which case
\begin{equation}
	(Xu)(t) := \sum_{j=1}^{n} u_j f_j(t). \label{eq:lin_comb_f}
\end{equation}
Then sampling $Xu(t)$ at $\{t_1,\dots,t_m\} \subset \mathcal{T}$, $t_1 < t_2 < \dots < t_m$, yields $y_i = Xu(t_i)$ and $v_{ij} = f_j(t_i)$
such that $Vu = y$. We will see, now, that the requirements of  $V^\transp, \Delta (V^\transp) \in \text{SC}_{m}$ in \cref{thm:unimodal} can be inherited from analogous properties of $X$.

In case that $\mathcal{T}$ contains finitely many elements $N$, we can assume that $X \in \mathds{R}^{N \times n}$ and all $m$-minors of $V$ and $\Delta (V^\transp)$ are, indeed, contained in $X^\transp_{[m]}$ and $[\Delta (X^\transp)]_{[m]}$. If $\mathcal{T}$ contains infinitely many elements, we can extend the notion of $X^\transp, \Delta X^\transp \in \sgc{m}$ by requiring that 
\begin{equation}
	X^\transp(t_1,\dots,t_m):= \begin{bmatrix}
		f_1(t_1) & f_1(t_2) & \dots & f_1(t_{m}) \\
		f_2(t_1) & f_2(t_2) & \dots & f_2(t_{m}) \\
		\vdots & \vdots & & \vdots \\
		f_n(t_1) & f_n(t_2) & \dots & f_n(t_m)
	\end{bmatrix}, \; \Delta X^\transp(t_1,\dots,t_m) \in \sgc{m} \label{eq:X_cond_for_thm}
\end{equation}
for all $t_1 < t_2 < \dots < t_m$, $t_i \in \mathcal{T}$, which trivially makes $V$ inherit these properties. In case that $\mathcal{T}$ is an interval, this may be tractable under the following more restrictive assumptions:
\begin{prop}\label[proposition]{prop:lin_comb_sc_m}
	Let $f_1,\dots,f_n \in C^{m-1}(a,b)$, $m < n$, $X: \mathds{R}^n \to C^{m-1}(a,b)$ be defined as in \cref{eq:lin_comb_f} and 
	\begin{equation}
		X_{\ast}^\transp (t) := \begin{bmatrix}
			f_1(t) & f_1^{(1)}(t) & \dots & f_1^{(m-1)}(t) \\
			f_2(t) & f_2^{(1)}(t) & \dots & f_1^{(m-1)}(t) \\
			\vdots & \vdots & & \vdots \\
			f_n(t) & f_n^{(1)}(t) & \dots & f_n^{(m-1)}(t)
		\end{bmatrix}.
	\end{equation}
Then, $X^\transp(t_1,\dots,t_m) \in \text{(S)SC}_m$ for all $a \leq t_1 < t_2 < \dots < t_m \leq b$, if
	\begin{enumerate}[label=\roman*.]
		\item for each fixed $k \in (1:m-1)$, the consecutive $k$-minors of $X_\ast^\transp (t)_{(:,(1:k))}$ share the same strict sign across all $t \in (a,b)$,
		\item the consecutive $m$-minors of $X_\ast^\transp (t)$ share the same (strict) sign for all $t \in (a,b)$.
	\end{enumerate}
	The same holds true, if we replace $X_\ast^\transp$ and $X^\transp$ with $\Delta X_\ast^\transp$ and $\Delta (X^\transp)$.
	\end{prop}
	
To illustrate \cref{prop:lin_comb_sc_m}, consider the following so-called \emph{confluent Vandermonde matrix} {\cite[Chapter~11]{respondek2025fast}} example. 
	\begin{exm}
		Let $f_j(t) := t^{j-1}$, $n = 4$ and $m = 3$. Then,
		\begin{equation*}
			X_\ast^\transp (t) = \begin{bmatrix}
				1& 0 & 0\\
				t & 1 & 0 \\
				t^2 & 2 t & 2\\
				t^3 & 3t^2 & 6t 
			\end{bmatrix} \quad \text{and} \quad \Delta X_\ast^\transp (t) = \begin{bmatrix}
				t-1 & 1 & 0\\
				t(t-1) & 2t-1 & 2 \\
				t^2(t-1) & t(3t-2) & 6t-2\\
			\end{bmatrix}
		\end{equation*}
		with
		\begin{align*}
			\det(X_\ast^\transp(t)_{((1:2),(1:2))}) &= 1, \; \det(X_\ast^\transp(t)_{((2:3),(1:2))}) = t^2, \; \det(X_\ast^\transp(t)_{((2:3),(1:2))}) = t^4\; \\
			\det(X_\ast^\transp(t)_{((1:3),(1:3))}) &= 2, \; \det(X_\ast^\transp(t)_{((1:3),(1:3))}) = 2t^3
		\end{align*}
		and 
		\begin{align*}
			\det(\Delta X_\ast^\transp(t)_{((1:2),(1:2))}) &= (t-1)^2, \; \det(\Delta X_\ast^\transp(t)_{((2:3),(1:2))}) = t^2(t-1)^2, \; \\
			\det(\Delta X_\ast^\transp(t)_{((1:3),(1:3))}) &= 2(t-1)^3.
		\end{align*}
		Using \cref{prop:lin_comb_sc_m}, this verifies that $$X^\transp(t_1,\dots,t_m), \; \Delta X^\transp(t_1,\dots,t_m) \in \text{SSC}_m.$$ However, by continuity of the determinant, it must follow then that $$X^\transp(t_1,\dots,t_m), \; \Delta X^\transp(t_1,\dots,t_m) \in \sgc{m}$$ if $\mathcal{T} = \mathds{R}_{\geq 0}$. As these derivations are also true for arbitrary $m$ and $n$, we can conclude that  \cref{thm:unimodal} is indeed a generalization of \cref{thm:log_conv}. 
	\end{exm}
	Similarly, it is possible to verify that also the family of real exponentials $f_i(t) = e^{\lambda_i t} $ with distinct $\lambda_i \in \mathds{R}$ on $\mathds{R}$ define such as an $X$. Since these families fulfill these properties for all $m \in (1:n)$ they are independent of the number of measurements available. 
	
	Finally, it should be noted that function families with $X^\transp(t_1,\dots,t_m) \in \sgc{m}$ are frequently found in the context of geometric modeling, e.g., B\'ezier curve, B-splines, etc. \cite{karlin1968total,hagen1992curve}. However, in this case $\Delta X^\transp(t_1,\dots,t_m) \in \sgc{m}$ may not be fulfilled for all $t_i$, despite an empirically observed unimodality of $p$, as the following example illustrates:
	\begin{exm}\label[example]{example:bernstein}
		Let $f_j(t) = \binom{10}{i-1} t^{i-1} (1-t)^{10-i+1}$, $j \in (1:11)$, 
		be the so-called \emph{Bernstein polynomials of order $10$}, $t_i = 0.1 i$, $i \in (1:4)$ and $V \in \mathds{R}^{4 \times 11}$ have entries $v_{ij} = f_j(t_i)$. Then, $p := (\normA{V_{(:,(1:4))}^{-1 }V_{(:,\{k\})}}{\ell_1})_{k \in (1:11)}\in \mathds{R}^{11}$ is unimodal (see~\cref{fig:p_bernstein}), despite 
		\begin{equation*}
			\det(\Delta(V^\transp)_{((1:4),:)} = 9.3076 \cdot 10^{-5} \quad \text{and} \quad \det(\Delta(V^\transp)_{(\{1,2,3,7\},:)} = -1.9686 \cdot 10^{-5},
		\end{equation*}
		i.e., $\Delta(V^\transp) \not \in \sgc{4}$. Moreover, $p$ drops below $1$ at $k = 10$, which verifies that our discussed failure phenomena is also present in a context outside of system matrices. In future work, we would like to close this gap to \cref{thm:unimodal} by developing less restrictive notions of variation bounding. 
	\end{exm}
	
	\begin{figure}[t]
		\centering
		\begin{tikzpicture}
			\begin{groupplot}[
				group style={
					group size=1 by 2,
					vertical sep=1.25 cm,
				},
				height=5 cm,
				width=12cm,
				]
				
				\nextgroupplot[xmin=1, xmax=11, 
				xlabel={$k$}, 
				ylabel={$p_k$},
				grid,
				grid,
				ytick={0.01,0.1,1,10},
				ymode=log
				]
				\addplot[line width=1.5pt, blue] 
				table [col sep=space,y expr=\thisrowno{1}]{p_bern_unimod.txt}; 
				
			\end{groupplot}
		\end{tikzpicture}
		
\caption{{Illustration of the vector $p := (\normA{V_{(:,1:m)}^{\dagger }V_{(:,\{k\})}}{\ell_1})_{k \in (1:n)}\in \mathds{R}^{11}$ 
	in \cref{example:bernstein}, where $V$ is constructed from the Bernstein basis polynomials of order $10$. 
	Despite the fact that  $\Delta (V^\transp) \not \in \sgc{4}$, the vector $p$ is still unimodal, indicating that unimodality can arise under less restrictive assumptions than those in \cref{thm:unimodal}.  Moreover, $p$ drops below $1$ at $k = 10$, showing that the failure phenomena discussed earlier can also occur in contexts beyond system matrices.}
}
		\label{fig:p_bernstein}
	\end{figure}
	\begin{rem}
	Note that \cref{thm:unimodal} can also be applied after a possible diagonal transformation $T = \diag(t_1,\dots,t_n)$ of $u$, $t_i \neq 0$, which results in replacing $V$ with $VT$. This is, because $\|T^{-1}u\|_{\ell_0} = \|u\|_{\ell_0}$, i.e., the transformation does not change \cref{eq:ell_0}. For example, if $A = \diag(x_1,\dots,x_m)$ with $x_i \leq 0$, and $b = \mathbf{1}_m$, then $V = \mathcal{C}^N(A,b) \not \in \sgc{m}$, but $VT \in \sgc{m}$ with $T = \diag(1,-1,1,\dots,(-1)^{n-1})$.
	\end{rem}

\section{Example -- Fuel Optimal Control}\label{sec:exm}
{In this section, we revisit the fuel-optimal regulator problem for a linear time-invariant system with fixed response time \cite{athans1966optimal}, as introduced in \cref{eq:ctrl_intro}.
Note that the original continuous-time problem formulation in \cite{athans1966optimal} also includes an upper bound \( |u(t)| \leq 1 \), which we chose to omit here.} %

As before, we can reformulate \cref{eq:ctrl_intro} into the form of \cref{eq:ell_0}, 
with
\begin{equation*}
	V = \mathcal{C}^N(A,b) \quad \text{and} \quad y = -A^{N}x(0).
\end{equation*}
and, after a possible state-space transformation, we may assume that $A = \diag(v)$, $v\in \mathds{R}^m$, and $b=\mathbf{1}_m$. To apply \cref{thm:log_conv,thm:uniquness}, we then require that $v > 0$: this condition guarantees that $V = \mathcal{C}^N(A,b) \in \text{STP}_m$, in which case $u^\ast$ is a unique solution to \cref{eq:ctrl_intro} if $\|u^\ast\|_{\ell_0} \leq \frac{m}{2}$ and, by our guaranteed unimodality property, one can efficiently determine whether there exists a $t^\ast \in (0:N-1)$ such that recovery of $u^\ast$ by \cref{eq:ell_1} necessitates $u^\ast(t) = 0$ for all $t \geq t^\ast$. To demonstrate this property, let us consider the following two concrete instances:
\begin{enumerate}[label=(\Roman*)]
	\item $A = \diag(\begin{bmatrix}
		0.98 & 0.97 & 0.96 & 0.95& 0.94
	\end{bmatrix}) $. \label{exm:mat_v1}
	{\item $A = \diag(\begin{bmatrix}
			0.8 & 0.7 & 0.6 & 0.5 & 0.4
		\end{bmatrix})
		$. \label{exm:mat_v2}}
\end{enumerate}
In \cref{fig:unimodal_vector_v1_v2,fig:unimodal_intro} and we can see that 
\begin{equation*}
	p_k =\normA{V_{(:,(1:m))}^{-1}V_{(:,\{k\})}}{\ell_1}, \; k \in (1:N)
\end{equation*}
is unimodal as predicted for $N = 500$ and $N =50$, respectively.  
\begin{figure}[t]
	\centering
	\begin{tikzpicture}
		\begin{groupplot}[
			group style={
				group size=1 by 1,
				vertical sep=1.25 cm,
			},
			height=5 cm,
			width=12cm,
			]
			
			\nextgroupplot[
			xmin=1, xmax=500,
			ymode=log,
			xlabel={$k$},
			ylabel={$p_k$},
			grid=both,
			minor ytick={},
			ytick={1,100,10000,1000000,10000000},
			yticklabel style={/pgf/number format/fixed},
			]
			
			\addplot[line width=1.5pt, blue]
			table[col sep=space, y expr=\thisrowno{1}]{unimodal_vector_v1.txt};
			
		\end{groupplot}
	\end{tikzpicture}
	
	\caption{Illustration of the predicted unimodality of the vector $p_k := \normA{V_{((:,(1:m))}^{-1}V_{(:,\{k\})}}{\ell_1}$ when $A$ is as in \labelcref{exm:mat_v1} and $k \in (1:500)$. $p_k$ remains above $1$ due to the narrowly clustered eigenvalues of $A$.}
	\label{fig:unimodal_vector_v1_v2}
\end{figure}
By comparison of both examples, it is notable that in the case of \labelcref{exm:mat_v1}, where the eigenvalues are narrowly clustered around $1$, $p_k$ never drops below $1$. Thus, our failure guarantees are inconclusive in this case. This observation is unsurprising as $p_{(m+1:N)}$ corresponds to an impulse response with very slow poles (see~\cref{thm:log_conv}). {Moreover, this behavior appears to be consistent with the continuous-time (CT) fuel-optimal control results, where controllability alone suffices to ensure the success of basis pursuit \cite{athans1966optimal,Nagahara2016handsoff}. In fact, one may interpret our system as a discretization of a CT LTI system with $A = e^{A_c h}$ for small step-size $h > 0$ and a Hurwitz stable diagonal $A_c \in \mathds{R}^{m \times m}$.} In contrast, for \labelcref{exm:mat_v2}, the sequence $p_k$ drops below $1$ at $k=36$ due to the presence of fast poles. According to \cref{thm:l1_necessary}, once the unimodal vector $p_k$ drops below $1$, it becomes impossible to recover an optimal solution $u^\ast$ to \cref{eq:ctrl_intro}, if any nonzero component $u^\ast(t)$ exists for $t \geq 36$. {Thus, revisiting our earlier observations in \cref{fig:ctrl_results_intro} from \Cref{sec:intro}: indeed, $u^\ast$ is the unique solution to \cref{eq:ctrl_intro} by \cref{thm:uniquness}; however \cref{thm:l1_necessary} and \cref{fig:unimodal_intro} guarantee that \cref{eq:ell_1}fails in recovering $u^\ast$.} On the other hand, if we use $A$ as in \labelcref{exm:mat_v1}, the optimal solution is numerically recovered when using a sufficiently high thresholding tolerance (namely $5 \cdot 10^{-8}$), which confirms our earlier intuition. {These observations suggest that a careful a priori selection of the horizon $N$ or the strategic placement of eigenvalues in $A$ may improve the success of basis pursuit in solving sparse optimal control problems. Note that this improvement occurs despite that standard conditions such as the restricted isometry property \cite{candes2005decoding} or mutual coherence \cite{donoho2003optimally} are not satisfied, due to the high cosine similarity between columns in $V$.} %

\section{Conclusion} \label{sec:cncl}
We have introduced a framework that provides tractable, deterministic failure guarantees for basis pursuit when solving sparse optimization problems under structural matrix constraints. Our analysis bridges necessary conditions from the dual approximation theorem with concepts from total positivity and control theory. {In particular, we focus on matrix structures -- including extended controllability, Hankel and Page matrices -- that frequently occur in control problems, as well as basis function from geometric modeling. The presented results reveal that the success of basis pursuit critically depends on the location of non-zero entries in the sparsest solution, or equivalently, on the placement of poles in system matrices.} Examples drawn from discrete-time optimal control problems further illustrate how these structural properties can influence performance in fuel-optimal control.

Overall, our work advances the theoretical understanding of sparse optimization and opens new avenues for practical design. We plan to extend our framework to other sparsity-inducing norms, such as the $\ell_2$-inducing sparsity norm considered in \cite{grussler2016lowrank,grussler2018low}, to quantify the gap between success and failure guarantees, and to broaden our results to other matrix structures. Finally, our findings on the log‐concavity of sequences of symmetric polynomials are also of independent algebraic interest.

\printbibliography

@Misc{amsmath,
  author =	 {{American Mathematical Society}},
  title =	 {User's Guide for the \texttt{amsmath} Package
                  (Version 2.0)},
  url =		 {ftp://ftp.ams.org/pub/tex/doc/amsmath/amsldoc.pdf},
  urldate =	 {2015-07-30},
  year =	 2002}

@Misc{pgfplots,
  author =	 {Christian Feuers\"anger},
  title =	 {Manual for Package \texttt{PGFPLOTS}},
  month =	 may,
  year =	 2015,
  url =		 {http://sourceforge.net/projects/pgfplots}
}

@preamble{"\newcommand{\noopsort}[1]{}"}

@string{IEEE-TAC = "IEEE Trans.\ Automat.\ Contr."}

@string{TAMS = "Trans.\ Amer.\ Math.\ Soc."}

@string{CUP  = "Cambridge University Press"}

@string{JWS  = "John Wiley {\&} Sons"}

@string{SIAM = "SIAM"}

@string{IEEE-CDC = "IEEE Conf. Decis. Control"}

@string{SIAM-MAA = "SIAM J.\ Matrix Anal.\ \& Appl."}

@string{IEEE-TIT = "IEEE Trans.\ Inform.\ Theory"}

@string{ACC = "Am. Control Conf."}

@string{IEEE-ICASSP = "Int. Conf. Acoust. Spee."}

@string{IEEE-TSPM = "IEEE Signal Proces. Mag."}

@string{IEEE-TSP = "IEEE Trans. Signal Proces."}

@article{sagan1992log,
  title={Log concave sequences of symmetric functions and analogs of the Jacobi-Trudi determinants},
  author={Sagan, Bruce E},
  journal=TAMS,
  volume={329},
  number={2},
  pages={795--811},
  year={1992}
}

@book{luenberger1968optimization,
	title={Optimization by Vector Space Methods},
	author={Luenberger, David G},
	year={1968},
	publisher=JWS,
}

@book{horn2012matrix,
title = "Matrix Analysis",
author = "Roger A. Horn and Charles R. Johnson",
year = "2012",
edition = "2",
publisher=CUP,
}

@article{grussler2016lowrank,
author = {Grussler, Christian and Giselsson, Pontus},
title = {Low-Rank Inducing Norms with Optimality Interpretations},
journal = {SIAM Journal on Optimization},
volume = {28},
number = {4},
pages = {3057-3078},
year = {2018},
	
}

@book{page,
  author = {A.A.H. Damen and P.M.J. Van den Hof and A.K. Hajdasinski},
  publisher = {Eindhoven University of Technology},
  title = {The Page Matrix: an excellent tool for noise filtering of Markov parameters, order testing and realization},
  year = {1982}
}

@book{antsaklis2006linear,
  title={Linear systems},
  author={Antsaklis, Panos J and Michel, Anthony N},
  year={2006},
  publisher={Springer}
}

@book{wedderburn1934lectures,
  title={Lectures on matrices},
  author={Wedderburn, Joseph Henry Maclagan},
  volume={17},
  year={1934},
  publisher={American Mathematical Soc.}
}

@article{pena,
  author          = {J.M. Peña},
  journal         = SIAM-MAA,
  number          = {4},
  title           = {Matrices with sign consistency of a given order},
  volume          = {16},
  pages           ={1100--1106},
  year            = {1995} 
}

@misc{roth2024system,
      title={On System Operators with Variation Bounding Properties}, 
      author={Chaim Roth and Christian Grussler},
      year={2024},
      eprint={2409.20275},
      archivePrefix={arXiv},
      primaryClass={math.OC},
}

@article{donoho2006compressed,
  title={Compressed sensing},
  author={Donoho, David L},
  journal=IEEE-TIT,
  volume={52},
  number={4},
  pages={1289--1306},
  year={2006},
}

@book{sra2011optimization,
  title={Optimization for machine learning},
  author={Sra, Suvrit and Nowozin, Sebastian and Wright, Stephen J},
  year={2011},
  publisher={MIT press}
}

@inproceedings{fardad2011sparsity,
  title={Sparsity-promoting optimal control for a class of distributed systems},
  author={Fardad, Makan and Lin, Fu and Jovanovi{\'c}, Mihailo R},
  booktitle="Proc.\ 2011 "#ACC,
  pages={2050--2055},
  year={2011},
  organization={IEEE}
}

@book{elad2010sparse,
  title={Sparse and redundant representations: from theory to applications in signal and image processing},
  author={Elad, Michael},
  year={2010},
  publisher={Springer Science \& Business Media}
}

@article{tillmann2013computational,
  title={The computational complexity of the restricted isometry property, the nullspace property, and related concepts in compressed sensing},
  author={Tillmann, Andreas M and Pfetsch, Marc E},
  journal=IEEE-TIT,
  volume={60},
  number={2},
  pages={1248--1259},
  year={2013},
  publisher={IEEE}
}

@inproceedings{chen2009sparse,
  title={Sparse LMS for system identification},
  author={Chen, Yilun and Gu, Yuantao and Hero, Alfred O},
  booktitle="Proc.\ 2009 "#IEEE-ICASSP,
  pages={3125--3128},
  year={2009},
}

@inproceedings{chen1994basis,
  title={Basis pursuit},
  author={Chen, Shaobing and Donoho, David},
  booktitle="Proc.\  28th Asilomar Conf.\  Signals, Systems and Computers",
  volume={1},
  pages={41--44},
  year={1994},
  organization={IEEE}
}

@book{karlin1968total,
	title={Total positivity},
	author={Karlin, Samuel},
	volume={1},
	year={1968},
	publisher="Stanford Univ.\ Press",
}

@book{foucart2013invitation,
  title={An invitation to compressive sensing},
  author={Foucart, Simon and Rauhut, Holger},
  year={2013},
  publisher={Springer}
}

@article{candes2005decoding,
  title={Decoding by linear programming},
  author={Cand{\`e}s, Emmanuel J and Tao, Terence},
  journal=IEEE-TIT,
  volume={51},
  number={12},
  pages={4203--4215},
  year={2005},
  publisher={IEEE}
}

@article{chandrasekaran2012convex,
  title={The convex geometry of linear inverse problems},
  author={Chandrasekaran, Venkat and Recht, Benjamin and Parrilo, Pablo A and Willsky, Alan S},
  journal={Foundations of Computational Mathematics},
  volume={12},
  number={6},
  pages={805--849},
  year={2012},
  publisher={Springer}
}

@article{donoho2003optimally,
  title={Optimally sparse representation in general (nonorthogonal) dictionaries via $\ell_1$  minimization},
  author={Donoho, David L and Elad, Michael},
  journal="Proc.\ National Academy of Sciences",
  volume={100},
  number={5},
  pages={2197--2202},
  year={2003},
  publisher={National Acad Sciences}
}

@article{candes2008introduction,
  title={An introduction to compressive sampling},
  author={Cand{\`e}s, Emmanuel J and Wakin, Michael B},
  journal=IEEE-TSPM,
  volume={25},
  number={2},
  pages={21--30},
  year={2008},
  publisher={IEEE}
}

@article{rauhut2010compressive,
  title={Compressive sensing and structured random matrices},
  author={Rauhut, Holger},
  journal={Theor. Found. Numer. Methods Spars. Recov.},
  volume={9},
  number={1},
  pages={92},
  year={2010}
}

@book{athans1966optimal,
  title={Optimal control: an introduction to the theory and its applications},
  author={M. Athans and P. L. Falb},
  year={1966},
  publisher={Mineola, NY: Dover Publications}
}

@ARTICLE{Nagahara2016handsoff,
  author={Nagahara, Masaaki and Quevedo, Daniel E. and Nešić, Dragan},
  journal=IEEE-TAC, 
  title={Maximum Hands-Off Control: A Paradigm of Control Effort Minimization}, 
  year={2016},
  volume={61},
  number={3},
  pages={735-747},
  }

@inproceedings{coulson2019regularized,
  title={Regularized and distributionally robust data-enabled predictive control},
  author={Coulson, Jeremy and Lygeros, John and D{\"o}rfler, Florian},
  booktitle="Proc.\ 58th "#IEEE-CDC,
  pages={2696--2701},
  year={2019},
}

@article{grussler2020variation,
title = {Variation diminishing linear time-invariant systems},
journal = {Automatica},
volume = {136},
pages = {109985},
year = {2022},
issn = {0005-1098},
author = {Christian Grussler and Rodolphe Sepulchre},
}

@article{grussler2021internally,
author = {Grussler, Christian and Burghi, Thiago and Sojoudi, Somayeh},
title = {Internally {H}ankel $k$-Positive Systems},
journal = {SIAM J. Control. Optim.},
volume = {60},
number = {4},
pages = {2373-2392},
year = {2022},
}

@article{gorodnitsky1997sparse,
  title={Sparse signal reconstruction from limited data using FOCUSS: A re-weighted minimum norm algorithm},
  author={Gorodnitsky, Irina F and Rao, Bhaskar D},
  journal=IEEE-TSP,
  volume={45},
  number={3},
  pages={600--616},
  year={1997},
  publisher={IEEE}
}

@book{vidyasagar2019introduction,
author = {Vidyasagar, M.},
title = {An Introduction to Compressed Sensing},
publisher = {Society for Industrial and Applied Mathematics},
year = {2019},
address = {Philadelphia, PA},
}

@article{grussler2018low,
  title={Low-rank optimization with convex constraints},
  author={Grussler, Christian and Rantzer, Anders and Giselsson, Pontus},
  journal=IEEE-TAC,
  volume={63},
  number={11},
  pages={4000--4007},
  year={2018},
}

@book{brunton2022data,
	title={Data-driven science and engineering: Machine learning, dynamical systems, and control},
	author={Brunton, Steven L and Kutz, J Nathan},
	year={2022},
	publisher={Cambridge University Press}
}

@book{hagen1992curve,
	title={Curve and surface design},
	author={Hagen, Hans},
	year={1992},
	publisher={SIAM}
}

@ARTICLE{zarmohdhigeojovTAC20,
	AUTHOR = {A. Zare and H. Mohammadi and N. K. Dhingra and T. T. Georgiou and M. R. Jovanovi\'c},
	TITLE = {Proximal algorithms for large-scale statistical modeling and sensor/actuator selection},
	JOURNAL = {IEEE Trans. Automat. Control},
	VOLUME = {65},
	NUMBER = {8},
	PAGES = {3441-3456},
	MONTH = {August},
	YEAR = {2020}
}

@INPROCEEDINGS{zarjovCDC18,
	AUTHOR = {A. Zare and M. R. Jovanovi\'c},
	TITLE = {Optimal sensor selection via proximal optimization algorithms},
	BOOKTITLE = {Proceedings of the 57th IEEE Conference on Decision and Control},
	YEAR = {2018},
	PAGES = {6514-6519}
}

@book{rugh1996linear,
author = {Rugh, Wilson J.},
title = {Linear system theory (2nd ed.)},
year = {1996},
isbn = {0134412052},
publisher = {Prentice-Hall, Inc.},
address = {USA}
}

@book{kwakernaak1991modern,
author = {Kwakernaak, Huibert and Sivan, Raphael and Strijbos, Rens C. W.},
title = {Modern signals and systems},
year = {1991},
isbn = {0138092524},
publisher = {Prentice-Hall, Inc.},
address = {USA}
}

@misc{carter2021complexity,
      title={The Complexity of Checking Partial Total Positivity}, 
      author={Daniel Carter and Charles Johnson},
      year={2021},
      eprint={2103.08742},
      archivePrefix={arXiv},
}

@Book{aitken1939determinants,
  author =       "A. C. Aitken",
  title =        "Determinants and Matrices",
  publisher =    "Oliver and Boyd",
  address =      "Edinburgh, UK",
  year =         "1944",
  edition = {Third},
  series =       "University mathematical texts",
}

@book{fallat2011totally,
  title={Totally Nonnegative Matrices},
  author={Fallat, Shaun M and Johnson, Charles R},
  year={2011},
  publisher={Princeton University Press},
  address={Princeton, NJ},
}

@misc{barshalom2022compound,
      title={Compound matrices in systems and control theory: a tutorial}, 
      author={Eyal Bar-Shalom and Omri Dalin and Michael Margaliot},
      year={2022},
      eprint={2204.00676},
      archivePrefix={arXiv},
}

@book{pinkus2009totally, place={Cambridge}, series={Cambridge Tracts in Mathematics}, title={Totally Positive Matrices}, publisher={Cambridge University Press}, author={Pinkus, Allan}, year={2009}, collection={Cambridge Tracts in Mathematics}}

@misc{grussler2026discrete,
      title={Discrete-Time Periodic Monotonicity Preserving Systems}, 
      author={Christian Grussler},
      year={2026},
      eprint={2503.23520},
      archivePrefix={arXiv},
}

@book{cullis1925matrices,
  title={Matrices and Determinoids},
  author={Cullis, Cuthbert Edmund},
  volume={2},
  year={1918},
  publisher={University Press}
}

@book{respondek2025fast,
	title={Fast matrix multiplication with applications},
	author={Respondek, Jerzy S},
	year={2025},
	publisher={Springer}
}
\appendix
\section{Appendix}
\label{sec:app}
\subsection{Proof to \cref{thm:l1_necessary}}
\begin{proof}
Without loss of generality, let $\|u^\ast\|_{\ell_1} = 1$. By \cref{lem:opt_cond}, $u^\ast$ can be a solution to \cref{eq:ell_1} only if there exists a vector $\beta^\ast\in\mathds{R}^mk$ satisfying $$
\|V^\transp \beta^\ast \|_{\ell_\infty} = u^\transp V^\transp\beta ^\ast =\|u^\ast\|_{\ell_1}.$$
In particular, this requirement implies that for each index $i \in (1:n)$, we must have $|V_{(\{i\},:)}^\transp\beta^\ast|\leq 1$ with equality (i.e, $V_{(\{i\},:)}^\transp \beta^\ast=\text{sign}(u_i^\ast)$) whenever $u_{i}^\ast\neq 0$. By our assumption that $V = V_{(:,(1:r))} V_{(:,(1:r))}^\dagger V$, this requires a $\tilde{\beta} \in \text{Im}\left(V_{(:,(1:r))}^\transp\right) \subset \mathds{R}^r$ such that  
\begin{equation*}
	|\tilde{\beta}_{j}| \leq 1, \;  j \in (1:r) \quad \text{and} \quad |{\tilde{\beta}}^\transp V_{(:,(1:r))}^\dagger V_{(:,\{i\})} | \leq 1
\end{equation*}
with equality whenever $u_{i}^\ast\neq 0$. Consequently, if $\|V_{(:,(1:r))}^{\dagger} V_{(:,i^\ast)}\|_{\ell_1} < 1$, such a vector cannot exist. %
\end{proof}

\subsection{Proof to \cref{thm:char_poly}}
\begin{proof}
	We only prove the case when $r = m$, since the other cases follow analogously. To this end note that $V_{(:,(1:m))}^{-1}V = \mathcal{C}(\tilde{A},\tilde{b})$ with $\tilde{b} := V_{(:,(1:m))}^{-1}b = e_1$ and
\begin{equation*}
	\tilde{A} := V_{(:,(1:m))}^{-1} A V_{(:,(1:m))} = \begin{bmatrix}
		0 & 0& 0 & \dots & 0 & -\alpha_0\\
		1 & 0 & 0 &\dots & 0 & -\alpha_1 \\
		0 & 1 & 0 &\dots & 0 & -\alpha_2\\
		\vdots & \vdots & \vdots    & & \vdots & \vdots\\
		0 & 0 & 0 &  \dots & 1 & -\alpha_{n-1}
	\end{bmatrix}.
\end{equation*}	
Therefore, if $\sum_{i=0}^{n-1} |\alpha_i| < 1$, it holds for the column sum norm (see, e.g., \cite[Example 5.6.4.]{horn2012matrix}) 
\begin{equation*}
	\|A\|_1 := \max_{\|x\|_{\ell_1} \leq 1} \|Ax\|_{\ell_1} = \max_{\|x\|_{\ell_1} = 1} \|Ax\|_{\ell_1}\leq 1,
\end{equation*}
and, in particular, $\|\tilde{A}^i\|_{1} \leq 1$ for all $i \geq 0$. However, since $$\|\tilde{A}^{m} \tilde{b}\|_{\ell_1} = \|\tilde{A}_{(:,\{m\})}\|_{\ell_1} <1,$$ it follows then that 
\begin{equation*}
	\forall k \geq 1: \; \|V_{(:,(1:m))}^{-1}V_{(:,\{k+m\})}\|_{\ell_1} = \|\tilde{A}^{m-1+k} \tilde{b}\|_{\ell_1} < 1.
\end{equation*}
	 \end{proof}

\subsection{Proof to \Cref{lem:spectrum}}
\begin{proof}
Since by the controllability assumption it holds that $\rk(V) = n$, it follows by \cref{lem:vb_m_sc_m} that $\mathcal{O}^N(A^\transp,b^\transp) \in \vb{m-1}$ for any $N \geq m$. Hence, for any $c^\transp \in \mathds{R}^{m}$ it must hold that $g^N =  \mathcal{O}^N(A^\transp,b^\transp) c^\transp$ fulfills $\vari{g^N} \leq m-1$ with $g^N_i = cA^{i-1}b$.  In other words, $g^N_i$ consists of the first $N$ samples of the impulse response to $(A,b,c)$. However, if there existed a $\lambda_i(A) \not \in \mathds{R}_{\geq 0}$, then by the controllability of $(A,b)$, it would be possible to choose a $c$ such that $G(z)$ (after possible pole-zero cancellation) only has poles in $\lambda_i(A)$ and $\overline{\lambda_i(A)}$. The impulse response to such a system changes its sign infinitely many times, which for sufficiently large $N$ contradicts $\vari{g^N} \leq m-1$. 
\end{proof}

\subsection{Proof to \cref{lem:pseudo_scm}}
\begin{proof}
Since $X^\dagger = (X^\transp X)^{-1}X^\transp \in \mathds{R}^{m \times n}$ due to our rank assumption, it follows that for any $\mathcal{I} \in \mathcal{I}_{n,m}$
\begin{equation*}
	\det(X^\dagger_{(:,\mathcal{I})}) = \frac{\det(X_{(\mathcal{I},:)})}{\det(X^\transp X)}.
\end{equation*}
Since $\det(X^\transp X) > 0$, as the determinant of a positive definite matrix, our claim holds.
\end{proof}
\subsection{Proof to \cref{thm:log_conv}}
\begin{proof}
	\emph{The case of $\rk(V) = m$:} We begin by showing that for each $k\ge 1$ and for every $i\in (1:m)$, we have $g_{(i)}(k) = W_k^\transp e_i \ge 0$.
	Since $\bar{V} := \mathcal{C}^N(A,b)^\transp \in \text{SC}_m$ for all $N > m$, and $\det(\bar{V}_{(:,(1:m))})\neq 0$ (by the controllability of $(A,b)$), \cref{lem:pena} implies that $
	\bar{V}_{(:,(m+k:N))} \bar{V}_{(:,(1:m))}^{-1} K_m$
	is totally positive. In particular, we get that
	$$(\bar{V}_{(:,(m+k:N))} \bar{V}_{(:,(1:m))}^{-1} K_m)^\transp e_i \geq 0$$ so that $g_{(i)}(k) \geq 0$ for all $k\ge 1$ and $i\in (1:m)$.
	Next, we prove the log concavity of $g_{(i)}(k)$ by showing that $g_{(i)}(k+1)^2 \geq g_{(i)}(k) g_{(i)}(k+2)$
	for every $k \geq 1$ and fixed arbitrary $i \in (1:m)$. To see this, set $c:=b^\transp$ and define the similarity matrix $$\tilde{A} := \mathcal{O}^m(A^\transp,c) A^\transp (\mathcal{O}^m(A^\transp,c))^{-1}.$$ Then, by \cref{lem:comp_properies} and the fact that $K_{[2]}= K_{[2]}^{-1}$, we get that
	\begin{align*}
		{\det \begin{bmatrix}
				g_{(i)}(k) & g_{(i)}(k+1)\\
				g_{(i)}(k+1) & g_{(i)}(k+2)
			\end{bmatrix}}&= {\det \begin{bmatrix}
				e_{m}^\transp \tilde{A}^k e_{i} & e_{m}^\transp \tilde{A}^{k+1} e_{i} \\
				e_{m}^\transp \tilde{A}^{k+1} e_{i} & e_{m}^\transp \tilde{A}^{k+2} e_{i} 
		\end{bmatrix}} = {L_{[2]}K_{[2]}K_{[2]}R_{[2]}} 
	\end{align*}
	where 
	\begin{align*}
		R := \begin{bmatrix}
			e_{i} & \tilde{A} e_{i}
			\end{bmatrix} \quad L := \begin{bmatrix}
			e_{m}^\transp \tilde{A}^k \\
			e_{m}^\transp	\tilde{A}^{k+1} 
			\end{bmatrix}.
	\end{align*}
		In particular, since $LK$ is a submatrix of the strictly totally positive matrix $$\mathcal{O}^{k}(A^\transp,c)A^{\transp^m} (\mathcal{O}^m(A^\transp,c))^{-1} K,$$ it holds that $(LK)_{[2]} > 0$ and by \cref{lem:comp_properies} $L_{[2]}K_{[2]} > 0$.  
		Thus, to conclude the desired inequality it remains to show that $K_{[2]}R_{[2]} =(KR)_{[2]} \leq 0$. To demonstrate this, we begin by noticing that 
		\begin{equation*}
			\tilde{A} =  \begin{bmatrix}
				0 & 1& 0 & \dots & 0 & 0\\
				0 & 0 & 1 &\dots & 0 & 0 \\
				\vdots & \vdots & \vdots    & & \vdots & \vdots\\
				0 & 0 & 0 & & 0 & 1\\
				(-1)^{m-1}a_{0}& (-1)^{m-2}a_{1}  & (-1)^{m-3}a_{2} &\cdots &-a_{m-2} & a_{m-1}
			\end{bmatrix}.
		\end{equation*}
		with characteristic polynomial
		$$\det(sI-\tilde{A}) = s^m - a_{m-1}s^{m-1} - \ldots - (-1)^{m-2} a_{1} s - (-1)^{m-1}  a_{0}.$$
		By \cref{lem:spectrum}, $\sigma(\tilde{A}) = \sigma(A) \subset \mathds{R}_{\geq 0}$
		and consequently $a_{0},\dots,a_{m-1} \geq 0$. Next, we observe that 
		$$
		KR = \begin{bmatrix}
			(-1)^{m-i+1}e_{m-i+1} & K\tilde{A}e_{m-i+1}
			\end{bmatrix}.$$
		This matrix has two non-zero values when $e_{i} = e_{1} $ and three non-zero entries otherwise, which is why one can verify that for all $1 \leq l_1 < l_2 \leq m$
		\begin{equation*}
			\det(KR_{(\{l_1,l_2\},:)}) = \begin{cases}
				-a_{m-1},& i=1\\
				-a_{m-i+1},& i\neq1\text{ and }l_1=1\\
				-1,&\text{otherwise} 
			\end{cases}.
		\end{equation*}
		Since these determinants are nonnegative, we conclude that $(KR)_{[2]} \leq 0$. This completes the proof of $g_{(i)}$ being log-concave.
		
		Finally, we show that the sequence $p(k):= \sum_{i=1}^m g_{(i)}(k)$ is log-concave. Since $p(k) \geq 0$ for all $k \geq 1$, it suffices to verify that $p(k+1)^2 \geq p(k)p(k+2)$.
		The proof follows analogously to the previous argument, with the only modification being that we redefine $R:= \begin{bmatrix}
			K \mathbf{1}_m & \tilde{A} K \mathbf{1}_m
			\end{bmatrix}.$
		This yields that
		\begin{equation*}
			{		KR = \begin{bmatrix}
					(-1)^{m-1}\phantom{-.---} & (-1)^{m-1} & \cdots & (-1)^{m-1}\\(-1)^{m}\sum_{i=1}^m a_{m-i} &(-1)^{m}\phantom{-}&\cdots&(-1)^{m}\phantom{-}\\    
				\end{bmatrix}^\transp},
		\end{equation*}
		which is why 
		\begin{equation*}
			\det(KR_{(\{l_1,l_2\},:)}) = 	\begin{cases}
				1 & \text{if } l_1 = -1 - \sum_{i=1}^m a_i\\
				0 & \text{else}
			\end{cases}
		\end{equation*}
		for all $1 \leq l_1 < l_2 \leq m$. Therefore, $(KR)_{[2]} \leq 0$, which completes the proof for the case of $\rk(V)=m$.\\
		\emph{The case of $\rk(V) = r < m$:} By our assumption and \cref{lem:pseudo_scm}, it holds that ${V}^\dagger_{(:,(1:r))} \in \sgc{r}$, which is why by \cref{lem:comp_properies}, $${V}^\dagger_{(:,(1:r))} {V} =  \mathcal{C}(V_{(:,(1:r))}^\dagger A, V_{(:,(1:r))}^\dagger b) \in \sgc{r}.$$ However, this means that we can simply apply our previous case to the full row rank $\mathcal{C}(V_{(:,(1:r))}^\dagger A, V_{(:,(1:r))}^\dagger b)$ to conclude the proof. 
		   \end{proof}
	
	\subsection{Proof to \cref{cor:transform}}
	\begin{proof}
	The proof follows directly from \cref{thm:log_conv} and the fact that $$
	K_m^\transp V_{(:,(1:m))}^{-1} V_{(:,\{m+k\})} = K_m^\transp \mathcal{C}(A,b)_{(:,(1:m))}^{-1} \mathcal{C}(A,b)_{(:,\{m+k\})} \text{ for all } k\geq 1,$$
	i.e., the $W_k$ coincide with those in \cref{thm:log_conv}.
	\end{proof}

	\subsection{Proof to \cref{thm:unimodal}}
		\begin{proof}
		To prove that $p_k$ is unimodal, it is sufficient to show that $\vari{\Delta p} \leq 1$. Since by definition $p_{(1:m)} = \mathbf{1}_m$, we only need to verify that $\vari{\Delta p_{(m:n)}} \leq 1$. To this end, define $W := V^\transp$ and $Q := W W_{((1:m),:)}^{-1}K_m$. Since $W \in \text{SC}_m$, it follows from \cref{lem:pena} that $Q_{((m+1:n),:)}$ is totally positive. Moreover, since by assumption $\Delta W \in \text{VB}_m$, it follows that $$
		\vari{\Delta Q\mathbf{1}_m} = \vari{\Delta W W_{((1:m),:)}^{-1}K_m \mathbf{1}_m} \leq m-1.$$
		In particular, the first $m$ rows of $\Delta Q \mathbf{1}_m$ contribute $m-2$ sign changes as seen from 
		\begin{equation}
			(\Delta Q)_{((1:m),:)} \mathbf{1}_m = \Delta K_m \mathbf{1}_m = \begin{bmatrix}
				\phantom{-}\cdots \phantom{-}& \phantom{-}2\phantom{-} & -2\phantom{-} & \phantom{-}2\phantom{-}
			\end{bmatrix}^\transp \label{eq:m_2_sign}.
		\end{equation}
		Since $p_{(m+1:n)} = Q_{((m+1:n),:)} \mathbf{1}_m$, it immediately follows that $$\vari{\Delta Q_{((m:n),:)} \mathbf{1}_m} =  \vari{\Delta p_{(m:n)}} \leq 1,$$ as desired. Finally, by \cref{eq:m_2_sign}, the one sign change in $\Delta p_{(m:n)}$ can only be due to a negative entry in $\Delta p_{(m+1:n)}$, from which it follows that $p$ is unimodal. 
		\end{proof}
		
		\subsection{Proof to \cref{prop:lin_comb_sc_m}}
		
		\begin{proof}
		We only need to show the claim for the case of $X^\transp$ as the case of $\Delta X^\transp$ follows analogously by considering the functions $$f_2-f_1,f_3-f_2,\dots,f_{n} - f_{n-1} \in C^{m-1}(a,b).$$ In \cite[Theorem~2.2.3. \& Supplement~2.2.1.]{karlin1968total}, this claim was proven for the case of $m=n$. Thus, applying this result to any selection of $m$ functions from $f_1,\dots,f_n$ results in demanding that for all $t \in (a,b)$ it must hold that
		\begin{enumerate}[label=\roman*.]
			\item $X_\ast^\transp(t)_{(:,(1:k))} \in \text{SSC}_k$, $k \in (1:m-1)$
			\item $X_\ast^\transp(t) \in \text{(S)SC}_m$.
		\end{enumerate}
	Then, using \cref{prop:consecutive_karlin} to verify these items permits us to only consider consecutive minors of $X_\ast^\transp(t)_{(:,(1:k))}$, which proves our claim. 
		\end{proof}

\end{document}